\documentclass{amsart}
\usepackage{epsfig}
\usepackage{amscd}
\usepackage{amsmath, amssymb, amsopn, amsthm, amsfonts}
\usepackage{graphics}
\usepackage{color}
\newtheoremstyle{thm-mod}
        {20pt} 
        {10pt} 
        {\itshape} 
        {-3pt} 
        {\bf} 
        {} 
        {\newline} 
        {\thmname{ #1}\thmnumber{ #2}\thmnote{ (#3)}} 
\theoremstyle{thm-mod}
\newtheorem{thm}{THEOREM}
\newtheorem{prop}{PROPOSITION}
\newtheorem{lem}[prop]{LEMMA}

\newtheoremstyle{def-mod}
        {10pt} 
        {10pt} 
        {} 
        {-3pt} 
        {\bf} 
        {} 
        {\newline} 
        {\thmname{ #1}\thmnumber{ #2}\thmnote{ (#3)}} 
\theoremstyle{def-mod}

\newtheorem{example}[prop]{Example}
\newtheorem*{eqn-label}{}

\def\reals{{\mathbb R}}
\def\cx{{\mathbb C}}

\def\Im{\,\mathrm{Im}\,}

\def\11{{\rm 1~\hspace{-1.4ex}l} }
\def\R{\mathbb R}

\begin{document}

\author{Wilhelm Schlag}
\title{Spectral Theory and Nonlinear PDE: A survey}
\address{Department of Mathematics, University of Chicago, 5734 South University Ave., Chicago, IL 60637, USA
} \email{schlag@math.uchicago.edu}
\thanks{The  author was partially supported by the NSF grant
DMS-0300081 and a Sloan Fellowship. This article is a much expanded
version of some lectures which the author held at the MSRI in
Berkeley during August 2005. He wishes to thank Jeremy Marzuola for
taking and typing notes and the MSRI and the organizers of the
dispersive PDE meeting for the opportunity to deliver these
lectures.}

\maketitle

\bigskip
This paper is not intended as a sweeping survey on the topic of
spectral theory and PDE as this would exceed any reasonable bounds.
Rather, we intend to outline some of the recent work by the author
and Joachim Krieger on the topic of stable manifolds for unstable
PDEs. More specifically, we intend  to give some ideas of how to
prove a recent result by Krieger and the author for the quintic wave
equation
\begin{eqnarray*}
\label{NLW} \Box \psi - \psi^5 = 0 \ \text{in} \ \R^3.
\end{eqnarray*}
This is the focusing, $H^1$--critical nonlinear wave equation (NLW).
There exist stationary solutions $\phi(x,a)$ of this equation given
by:
\begin{eqnarray*}
\phi(x,a) = (3a)^{\frac{1}{4}} (1 + a|x|^2)^{-\frac{1}{2}}, \qquad
a>0
\end{eqnarray*}
The question is what happens if the initial data of NLW is a small
perturbation of $\phi$. To answer this question, we first linearize
around $\phi$. The ansatz $\psi(x,t) = \phi(x,a) + u(x,t)$ leads to
the equation
\begin{eqnarray*}
\partial_{tt} u - \Delta u + V u = N(u,\phi),
\end{eqnarray*}
with $V = -5 \phi^4$ and a nonlinear term \[ N(u,\phi)=10
u^2\phi^3+10u^3\phi^2 + 5u^4\phi + u^5
\]
Let $H = - \Delta + V$. It is known that the spectrum of such a
Schr\"odinger operator on the interval $(0,\infty)$ is purely
absolutely continuous. Also, $H$ has a simple negative eigenvalue,
say $-k^2$ (the ground-state). Indeed,
\[ \inf_{f\in H^1,\|f\|_2=1} \langle Hf,f\rangle \le \langle H\phi,\phi\rangle = -4\int_{\R^3} \phi^6\, dx <0\]
The presence of negative spectrum implies that the static
solutions~$\phi(\cdot,a)$ are linearly (exponentially) unstable. We
remark that by Agmon's estimate, we know that the ground-state
function $g$, which satisfies $H g = -k^2 g$,  decays exponentially.
In addition, $g > 0$, and $g$ is radial. Because of symmetries, $H$
has a zero mode. More precisely, by the dilation and translation
symmetries
\[ H \partial_a \phi(\cdot,a) =0, \qquad H \nabla \phi(\cdot,a) =0
\]
It is not an easy task to show that $H$ has exactly one negative
eigenvalue (which then must be simple), and no other zero modes than
those four which we have specified. Numerically, this has been
verified by means of the Birman-Schwinger method (the numerics show
that the Birman-Schwinger operator has exactly five eigenvalues
$\ge1$ counted with multiplicity. This corresponds precisely to the
four independent eigenfunctions and the resonance function).

However, if we restrict $H$ to the invariant subspace of radial
functions, then the determination of the spectrum of this restricted
operator is elementary. Indeed, $H f = Ef$ for radial $f=f(r)$ is
the same as $(-\frac{d^2}{dx^2} + V)rf(r)=E\,rf(r)$. In other words,
we are reduced to an eigenvalue problem for the half-line operator
$\tilde H=-\frac{d^2}{dx^2} + V$ with a Dirichlet condition at
$r=0$. In our case, the function $r\partial_r \phi(r,a)$ has a
unique positive zero and solves $\tilde H (r\partial_r \phi(r,a))=0$
(moreover, it is bounded and therefore a resonance function
for~$\tilde H$). By Sturm theory, it follows that there is a unique
negative eigenvalue which is simple. It agrees with $-k^2$ from
above. Hence, the operator~$\tilde H$ has a one-dimensional
exponential instability (in the radial setting) and the evolution of
the wave equation driven by~$H$ has a codimension one stable
manifold. Our main theorem states that this remains true (in some
sense) on the nonlinear level.

\begin{thm}
There exists a $\delta>0$ such that the following statement holds.
Let $B_\delta (0)$ be the ball centered at $0$ of size $\delta$ in
$H^3_{rad}(\R^3) \times H^2_{rad}(\R^3)$. There exists a
codimension~$1$ Lipschitz manifold $\Sigma$ which divides
$B_\delta(0)\setminus\Sigma$ into two connected components and such
that if $(u_1,u_2) \in \Sigma \cap B_{\delta} (0)$, then the Cauchy
problem for NLW with data
\[(\psi(\cdot,0),\partial_t \psi(\cdot,0))=(\phi(\cdot,1)+u_1,u_2)\]  has a global solution
$\psi(x,t)$. Moreover, there is the representation
$$\psi(x,t) = \phi(x,a(\infty)) + u(x,t)$$ where $|a(\infty)- 1|\lesssim \delta$,
and $\|u(\cdot,t)\|_\infty \lesssim \delta \langle t\rangle^{-1}$ as
well as
 \[(u,\partial_t u) = (v, \partial_t v) + o_{\dot H^1\times L^2}(1)
 \text{\ \ as\ \ }t\to\infty\]
  with $\Box {v} = 0$ and
$({v},\partial_t {v})(0) \in \dot{H}^1 \times L^2$.
\end{thm}

One motivation for this problem came from a numerics paper by Bizon,
Chmaj, and Tabor which appeared in {\em Nonlinearity} in~2004. These
are physicists who are mainly interested in the Einstein equations
of general relativity, but proposed the quintic NLW as a model case.
They predicted numerically that the stable manifold should exist, at
least in a small neighborhood of $\phi$. Moreover, Bizon et.\ al.\
made conjectures about the behavior above and below the manifold
(one should lead to blow-up, the other to dispersion).

Another motivation for studying the quintic NLW came from the
author's work on the focusing nonlinear Schr\"odinger
equation~(NLS), which we now describe. Consider the equation
\begin{eqnarray*}
i \partial_t \Psi + \Delta \Psi + | \Psi |^{2 \sigma} \Psi = 0 \
\text{in} \ \R^3,
\end{eqnarray*}
for $\sigma < 2$ (where $2$ is the $H^1$ critical exponent in
$\R^3$. We also have $\sigma_{crit} = \frac{2}{d} = \frac{2}{3}$,
which is the $L^2$ critical exponent). The Hamiltonian and the
charge are
conserved functionals: \\
\begin{eqnarray*}
\mathcal{H} ( \Psi (t) ) & = & \int \left( \frac{1}{2} | \nabla \Psi |^2 - \frac{1}{2(\sigma+1)} | \Psi |^{2(\sigma+1)} \right) dx \\
Q( \Psi(t)) & = & \frac{1}{2} \int | \Psi |^2 dx.
\end{eqnarray*}
For $\sigma < \frac{2}{d}$, we get global existence in $H^1$ using
the Gagliardo-Nirenberg inequality and the above conservation laws.
  For $\sigma \geq \frac{2}{d}$, we in fact can have blow-up in finite time.  This result is easy to see simply by looking at the
  so-called {\em virial identity}. This approach is known as Glassey's
  argument
  and is based on the functional \[\int |x|^2 | \Psi (t,x) |^2 dx.\]
  Indeed, if $\sigma\ge 2/d$, then
  \[ \frac{d^2}{dt^2} \int |x|^2 | \Psi (t,x) |^2 dx \le 16
  \mathcal{H}(\Psi(0))\] Thus,
   any solution of negative energy must blow up in finite time (at
   least in the sense of the $H^1$ norm). Indeed, at some point the
   positive quantity above disappears and the solution must blow-up.
   Strictly speaking, to
   apply the virial identity, we need to work in the space $\langle
   x\rangle^{-1}H^1$ (it suffices to assume that $\Psi(0) \in
\langle x \rangle^{-1} H^1$).

There exist stationary state solutions for NLS as well. Make the
ansatz $\Psi(t,x) = e^{i t \alpha^2} \phi$. Then $\phi$ satisfies
the equation
\begin{eqnarray*}
(\alpha^2 - \Delta) \phi = | \phi |^{2 \sigma} \phi.
\end{eqnarray*}
Berestycki and Lions showed using variational methods that there
exist solutions to such equations that are positive, radial, and
exponentially decaying.  In many cases, uniqueness of said
ground-state solutions has been proved by Coffman, Kwong, and
McLeod. $\phi$ must satisfy the Derrick--Pohozaev identities and for
monomial nonlinearities we have a scaling property, namely
$\phi(x,\alpha) = \alpha^{\frac{1}{\sigma}} \phi (\alpha x, 1).$

There is a group of natural symmetries acting on NLS. There is the
modulation symmetry $\Psi \to e^{i \gamma} \Psi$ and  the Galilei
transforms
$$\mathcal{G} (t) = e^{i(-v^2 t + x \cdot v )} e^{-ip \cdot (y + 2tv)},$$
for $p = -i \frac{\partial}{\partial_x}$. These are the unique
transformations that satisfy
\[ \mathcal{G}(t) e^{it\Delta} = e^{it\Delta} \mathcal{G}(0)\]
with $\mathcal{G}(0)=e^{ix \cdot v } e^{-ip \cdot y }$. Note that
this latter transform at time $t=0$ is the Galilei transform from
classical mechanics, i.e., $x\mapsto x+y$, $p\mapsto p+v$.
 Finally, we have dilation symmetries
from the scaling relation above for ground state solutions.  In
total, we have a $2d+2$ vector of symmetry parameters
$$ \pi = ( \gamma, v, y, \alpha ).$$

The relevance of these symmetries can be seen in the following
theorem of Michael Weinstein.
\begin{thm}
For any $\epsilon$, there exists a $\delta$ such that if
 $$\inf_{\gamma,y} \| \Psi(\cdot,0) - e^{i \gamma} \phi(\cdot + y,\alpha) \|_{H^1} < \delta,$$
then
$$\inf_{\gamma,y} \| \Psi(x,t) - e^{i \gamma} \phi( \cdot + y,\alpha) \|_{H^1} < \epsilon.$$
\end{thm}
It is essential to mod out the symmetries in the second inequality,
since a small perturbation of the ground-state $\phi$ (which
unperturbed has zero momentum) will impart a nonzero (but very
small) momentum on it. Consequently, the soliton will have moved a
large distance in finite (but large) time. In a similar fashion, the
phase will also change by size one in finite time.

Weinstein's theorem expresses what is known as {\em orbital
stability} of the ground-state solutions. During the late 1980s
Grillakis, Shatah, and Strauss developed a fundamental theory of
orbital stability of standing waves generated by symmetries. Their
theory applies to general Hamiltonian equations which are invariant
under the action of Lie groups. In particular, it gives the
dichotomy $\sigma<2/d$ and $\sigma\ge 2/d$ for the focusing NLS, but
applies to a much wider range of examples.

Note that the conjectured instability of excited state solutions is
largely open. An excited state is a solution of \[ -\Delta \phi
+\alpha^2 \phi= |\phi|^{2\sigma}\phi\] which changes sign. Important
work has been done in that direction by T.~P.~Tsai and H.~T.~Yau for
NLS with a potential. Generally speaking, NLS with a potential or
Hartree equations is a fascinating field to which we cannot do
justice here, see for example the work of Fr\"ohlich et.~al.

Berestycki and Cazenave showed that for $\sigma \geq \frac{2}{d}$,
blow-up can occur for arbitrarily small perturbations of ground
states. This is easy for $\sigma = \frac{2}{d}$ using Glassey's
argument from above and the fact that $\mathcal{H} ( (1+\epsilon)
\phi ) < 0$ for $\epsilon>0$ since $\mathcal{H} (\phi) = 0$.

Now, let us discuss the spectral analysis of linearized operators
for NLS.  Set $\psi = e^{i t \alpha^2} \phi + R$ and linearize. Note
that we are assuming $\alpha$ does not change with time here, but we
do not let that concern us at this point.  We end up with an
equation
\begin{eqnarray*}
i \partial_t \left( \begin{array}{l}
R \\
\bar{R}
\end{array} \right) + H \left( \begin{array}{l}
R \\
\bar{R}
\end{array} \right) = N,
\end{eqnarray*}
for $N$ a nonlinear term quadratic in $R$ and
\begin{eqnarray*}
H = \left( \begin{array}{cc}
\Delta - \alpha^2 + (\sigma + 1) \phi^{2 \sigma} & \sigma \phi^{2 \sigma} \\
- \sigma \phi^{2 \sigma} & - \Delta + \alpha^2 - (\sigma + 1)
\phi^{2 \sigma}
\end{array} \right).
\end{eqnarray*}
If $\phi = 0$, we have a self-adjoint operator with purely
absolutely continuous spectrum given by $(-\infty,-\alpha^2] \cup
[\alpha^2,\infty)$. If we do turn on $\phi$, then this set turns
into the essential spectrum of $H$ by Weyl's criterion. In addition,
there is an eigenvalue at $0$ due to the symmetries. Indeed, writing
NLS as $\mathcal{N} ( \Psi_\pi) = 0$ where $\Psi_\pi$ indicates the
symmetries given by the vector $\pi$ acting on $\Psi$, we conclude
that $\mathcal{N}'(\Psi_\pi)
\partial_\pi \Psi_\pi = 0$. Since clearly $H(\Psi_\pi) = \mathcal{N}' (\Psi_\pi)$ by definition of the
linearized operator, it follows that  $\partial_\pi \Psi_\pi$ is a
zero energy state of $H$. In the NLS case, they are (generalized)
eigenfunctions with zero eigenvalue (but for the quintic NLW from
above, there is a resonance and not just an eigenvalue at zero). If
we look at $R = u + iv$ and write $H$ in terms of $u$, $v$, then we
obtain the equivalent matrix operator
\begin{eqnarray*}
H = \left( \begin{array}{cc}
0 & i L_{-} \\
-i L_{+} & 0
\end{array} \right),
\end{eqnarray*}
where
$$L_{-} = - \Delta + \alpha^2 - \phi^{2 \sigma}$$
and
$$L_{+} = -\Delta + \alpha^2 - (2 \sigma + 1) \phi^{2 \sigma}.$$
Clearly, $L_{-} (\phi) = 0$ and $L_{+} (\nabla \phi) = 0$.  Hence,
we have
\begin{eqnarray*}
\left( \begin{array}{l}
0 \\
\phi \end{array} \right), \left( \begin{array}{l}
\nabla \phi \\
0 \end{array} \right) \in \ \text{Ker} \ (H).
\end{eqnarray*}
Differentiating $L_-\phi=0$ in $\alpha$ yields  $L_{+}
(\partial_\alpha \phi) = - 2 \alpha \phi$. Therefore,
\begin{eqnarray*}
H \left( \begin{array}{l}
\partial_\alpha \phi \\
0 \end{array} \right) = \left( \begin{array}{l}
0 \\
2i \alpha \phi
\end{array} \right).
\end{eqnarray*}
In fact, it is easy to see that
\begin{eqnarray*}
\left( \begin{array}{l}
0 \\
x \phi
\end{array} \right),
\left( \begin{array}{l}
\partial_\alpha \phi \\
0
\end{array} \right) \in \ \text{Ker} \ (H^2).
\end{eqnarray*}
This shows that the eigenvalue $0$ is  of geometric multiplicity at
least $d+1$ and of algebraic multiplicity at least $2d+2$, but we
actually have equality here (at least if $\sigma\ne\frac{2}{d}$) due
to a result of M.~Weinstein from~$1986$. So, $H$ is truly a non-self
adjoint operator with a nontrivial Jordan form at zero energy. In
the critical case $\sigma=\frac{2}{d}$ the dimension of the root
space at zero increases to $2d+4$ due to the extra pseudo-conformal
symmetry.

Let us now analyze the operator $e^{i t H}$.  Using the $L_{-}$,
$L_{+}$ version of $H$, we see that $H f = E f$ implies
\begin{eqnarray*}
H^2 \left( \begin{array}{l}
f_1 \\
f_2
\end{array} \right) = \left( \begin{array}{cc}
L_{-} L_{+} & 0 \\
0 & L_{+} L_{-}
\end{array} \right) \left( \begin{array}{l}
f_1 \\
f_2
\end{array} \right) = E^2 \left( \begin{array}{l}
f_1 \\
f_2
\end{array} \right).
\end{eqnarray*}
Since $L_{-} \phi = 0$ and $\phi > 0$ it follows that $L_{-}$ is a
non-negative operator with $\phi$ as a ground state. By work of
M.~Weinstein,  $L_{+}$ has exactly one negative eigenvalue and the
kernel is spanned by the derivatives $\partial_j\phi$.  From the
above argument, since $E^2 f_1$ is in the range of $L_{-}$, we know
that $f_1 \perp \phi$.  Therefore, we have
\begin{eqnarray*}
\sqrt{L_{-}} L_{+} \sqrt{L_{-}} \tilde{f_1} = E^2 \tilde{f_1},
\end{eqnarray*}
for $\tilde{f_1} = L_{-}^{- \frac{1}{2} } f_1$.  Thus, $E^2$ is real
as an eigenvalue of a self-adjoint operator. The dichotomy of linear
stability or instability is decided by the question: When is $E^2 <
0$? In order to answer this, look at
\begin{eqnarray*}
\inf \frac{\langle \sqrt{L_{-}} L_{+} \sqrt{L_{-}} f,f \rangle}{\langle f,f \rangle} & < & 0 \\
\Leftrightarrow \inf_{\tilde{f} \perp \phi, \| \tilde{f} \|_2 = 1}
\langle L_{+} \tilde{f},\tilde{f} \rangle & < & 0.
\end{eqnarray*}
Denote this latter infimum as $\mu_0$ and assume that $\mu_0<0$. One
can show that the minimum is attained at some function $f_0$. Using
Lagrange multipliers, we have \[L_{+} f_0 = \mu_0 f_0 + \lambda
\phi,\quad f_0 \perp \phi,\qquad \mu_0 = \langle L_{+} {f_0},{f_0}
\rangle.\]
 If $\lambda = 0$, then $\mu_0$ is
a negative eigenvalue, hence the ground state of $L_{+}$ and this
contradicts the orthogonality relation.  So, we have $f_0 = (L_{+} -
\mu_0)^{-1} \phi$.  Consider the function
\begin{eqnarray*}
h(\mu) = \langle (L_{+} - \mu)^{-1} \phi,\phi \rangle.
\end{eqnarray*}
Note that this goes to $-\infty$ for $\mu$ tending to the unique
negative eigenvalue of $L_+$ and $h(0)$ makes sense since the kernel
of $L_{+}$ is orthogonal to $\phi$ (remember it is $\nabla \phi$).
We know $h'(\mu) > 0$ as it is $\| ( L_{+} - \mu )^{-1} \phi
\|^2_2$. Hence,
 $h(0) > 0$ iff $\mu_0 < 0$.  But,
\begin{eqnarray*}
h(0) & = & -\frac{1}{2 \alpha} \langle \partial_{\alpha} \phi,\phi \rangle \\
& = & - \frac{1}{4 \alpha} \partial_\alpha \| \phi_\alpha \|^2_{L^2}
> 0,
\end{eqnarray*}
Finally, one checks from scaling that $\partial_\alpha \|
\Psi_{\alpha} \|^2 < 0 \iff \sigma
> \frac{2}{d}$. This (up to the equality sign) is the well-known instability condition
for NLS. We now state our first stable manifold theorem for
focusing, cubic NLS in three dimensions. It is conditional on the
assumption that there are no imbedded eigenvalues in the essential
spectrum of the linearized operator.

\begin{thm}
(Schlag - 2004) \\
Consider NLS with $\sigma=1$ in $\R^3$.  Assume that the linearized
operator does not have any imbedded eigenvalues in its essential
spectrum. Then there exists $\delta>0$ such that the following
holds: Let $B_\delta (0) \subset W^{1,1} \cap W^{1,2}$ be the ball
centered at $0$ of radius~$\delta$. There exists a codimension one
Lipschitz graph $\Sigma\subset B_\delta(0)$ which divides
$B_\delta(0)\setminus\Sigma$ into two connected components and so
that for initial data in $\Psi(\cdot,0)\in \Sigma+\phi(\cdot,1)$,
there exists a global solution
$$\Psi(x,t) = W(x,t) + R(x,t),$$ with $\| R(\cdot,t) \|_{\infty}
\lesssim \delta \langle t \rangle^{-\frac{3}{2}}$, $R = e^{i t
\Delta} \tilde{R} + o_{L^2} (1)$.  We have
$$W(x,t) = e^{i \theta(t,x)} \phi(x - y(t), \alpha(t)),$$
$$\theta(t,x) = -\int^t_0 v^2(s) dx + v(t) x + \gamma(t)$$ and
$$y(t) = 2 \int_0^t v(s) ds + y_0 (t).$$
Finally, the vector $\pi(t):=(\gamma(t), v(t), y_0(t), \alpha(t))$
approaches a terminal vector $\pi(\infty)$  which is~$\delta$ close
to the initial vector $\pi(0)=(0,0,0,1)$.
\end{thm}

It is clear that we cannot have stability for initial data in the
entire ball $\phi+B_\delta(0)$. In fact, $\sigma=1$ is in the
supercritical range for which Berestycki and Cazenave have exhibited
arbitrarily small perturbations of $\phi$ that lead to finite-time
blow-up. On the other hand, we do not understand the general
behavior of solutions with data in
$\phi+B_\delta(0)\setminus\Sigma$.

We remark that the method of proof does not extend to the entire
supercritical range $\sigma>2/3$. More precisely, an important gap
property for the operators $L_{\pm}$ (which we defined above) breaks
down below some value $\sigma_*=0.91396\ldots$. The gap property
here refers to  $L_+$ and $L_-$ having neither an eigenvalue in
$(0,\alpha^2]$ nor a resonance at the edge~$\alpha^2$.  This was
shown in some recent numerical work by Laurent Demanet and the
author. This numerical work also establishes that the gap property
does hold above $\sigma_*$, and thus in particular for $\sigma=1$
(which is essential for the proof of the theorem).

For the proof of the theorem, we follow the basic strategy of the
modulation theory of Soffer and Weinstein. More precisely, we seek
an ODE for $\pi(t) = (\gamma,y,v,\alpha)(t)$ and a PDE for $R$. The
ODE is constructed in such a way that the non-dispersive nature of
the root-space of $H$ is eliminated. Note that the dimension of the
root-space is $8$, as explained above, and $\pi(t)$ has eight
components. This agreement is of course not only essential but also
no coincidence, as $\pi$ has as many components as there are
symmetries, and the dimension of the root-space must also agree with
the number of these symmetries.

Of fundamental importance are the following estimates on the
linearized operator:
\begin{eqnarray*}
\| e^{itH} P_s f \|_\infty & \lesssim & t^{-\frac{3}{2}} \| f \|_1,  \\
\| e^{itH} P_s f \|_2 & \lesssim & \| f \|_2
\end{eqnarray*}
where $P_s$ projects onto the stable part of the spectrum (i.e., it
is the identity minus the Riesz projection onto the discrete
spectrum). The proof of these bounds not only requires that there
are no imbedded eigenvalues in the essential spectrum, but also the
following properties:

\begin{itemize}
\item zero is the only point in the discrete spectrum
\item the edges $\pm\alpha^2$ are neither eigenvalues nor resonances
\end{itemize}

We prove these facts by modifying some arguments of Galina Perelman
for the case $d=1$ and $\sigma=2$. It is at this point that the
aforementioned numerical work by Demanet and the author  becomes
essential.

In one dimension a corresponding result was proven by Joachim
Krieger and the author in the entire supercritical range $\sigma>2$.
This result was {\em completely unconditional} and also did not rely
on any numerics. Rather, one can prove all the necessary spectral
properties by means of the scattering theory established by Buslaev
and Perelman, as well as by invoking the explicit form of the ground
state in~$d=1$ (the latter is needed to prove that there are no
imbedded eigenvalues). The linear estimates needed for $d=1$ are

\begin{eqnarray*}
\|e^{itH} P_s f\|_\infty &\lesssim t^{-\frac12} \|f\|_1 \\
\| \langle x \rangle^{-1} e^{itH} P_s f \|_{\infty} & \lesssim
t^{-\frac{3}{2}} \| \langle x \rangle f \|_1
\end{eqnarray*}

The latter bound is somewhat curious, in so far as it does not hold
in the free case $H=-\partial^2$. Rather, it exploits that the edges
of the essential spectrum are not resonances (which of course fails
in the free case). Heuristically, one can arrive at the second bound
as follows. Suppose the scalar operator $-\partial^2 +V$ does not
have a resonance at zero energy. This is well-known to be equivalent
to the fact that the distorted Fourier basis
$\{e(x,\xi)\}_{\xi\in\R}$ vanishes when $\xi=0$. Now write the
evolution relative to this basis:
\begin{eqnarray*}
e^{itH}P_s f(x)= \int e^{it \xi^2} e(x,\xi) \langle f
,e(\cdot,\xi)\rangle \, d\xi.
\end{eqnarray*}
By the vanishing of $e(\cdot,0)$, we integrate by parts in $\xi$
gaining a factor of $t^{-1}$. However, this evidently comes at the
expense of a factor of~$x$. Hence, we obtain the $t^{-\frac32}$
decay by invoking the dispersive decay on top of the $t^{-1}$ gain,
leading to the above bound. On a more heuristic level, we remark
that compact data $f$ should be thought of as having traveled
distance $t$ after time $t$. Due to the weight $\langle
x\rangle^{-1}$ this will lead to a gain of $t^{-1}$ on top of the
dispersive decay of~$t^{-\frac12}$. The point of not having a
resonance at the edge is to avoid having some mass "being stuck"
around the origin. This would of course invalidate the heuristic
idea of transport.

The possibility of obtaining improved decay bounds at the expense of
weights was apparently discovered by Murata in the~$1980s$. It
applies only to those cases where the free operator has a zero
energy resonance, i.e., when $d=1,2$. In dimension two, Murata
discovered an improvement of $(\log t)^{-2}$ over the free
dispersive decay.

We remark that the presence of these weights leads to substantial
technical problems in the proof of the 1-d stable manifold theorem
(as compared to the 3-d result).

Finally, we would like to point the reader to the second stable
manifold paper with Krieger which deals with the pseudo conformal
blow-up solutions for critical NLS in $d=1$.  The $t^{-\frac32}$
bound again plays a crucial role, but this is not the venue to
explain the blow-up argument in any detail.

\bigskip As stated above, we desire to understand stability of
the focusing, quintic NLW
\begin{eqnarray*}
\Box \psi - \psi^5 = 0
\end{eqnarray*}
around the special solutions $\phi(\cdot,a)$.  Recall that  Aubin
showed that there exists $\phi = \phi(\cdot,a)>0$ such that $-
\Delta \phi (\cdot,a) - \phi^5 (\cdot,a) = 0$. Explicitly, with
$a>0$
\[ \phi(x,a)=(3a)^{\frac14}(1+a|x|^2)^{-\frac12} \]
 Note from the formula above
that $\phi(x,a) \sim \frac{1}{|x|}$ and also $\partial_a \phi(x,a)
\sim \frac{1}{|x|}$. In particular, neither $\phi$ nor $\partial_a
\phi$ belong to~$L^2(\R^3)$. We also remark that it is essential for
us to work with the quintic, i.e., $H^1$-critical equation for a
number of reasons. For example, by a result of Gidas and Spruck,
there are no positive solutions of $\Delta \phi + \phi^p=0$ in
$\R^3$ if $p<5$.

We wish to study radial perturbations to the initial data of the
form
\begin{eqnarray*}
\left( \begin{array}{l}
\psi(0) \\
\partial_t \psi(0)
\end{array} \right) = \left( \begin{array}{l}
\phi(\cdot,1) \\
0
\end{array} \right) + \left( \begin{array}{l}
f_1 \\
f_2
\end{array} \right).
\end{eqnarray*}
As mentioned before, Bizon, Chmaj, and Tabor, as well as Szpak,
studied this problem numerically. They tested initial data $(f_1,
f_2)$ from a one-parameter family of functions (e.g.,~Gaussians with
changing variance) and observed that there is a unique point on such
a curve of data at which there is long-time stability. This point
splits the curve into two connected components (i.e., intervals).
One interval led to instability (or possibly blow-up), whereas
another led to scattering to a free wave. By letting this curve of
initial data vary, they then observed that the unique stability
point  sweeps out a codimension one stable manifold on which the
family $\phi(\cdot,a)$ acts as a one-dimensional attractor.

Our theorem establishes the existence of such a codimension one
stable manifold. However, at this point we cannot describe the
behavior off the manifold. We now state the theorem, which already
appeared above, in full detail.

\begin{thm}
(Krieger-Schlag - 2005) \\
Let $R>1$.  Denote
$$X_R = \{ (f_1,f_2) \in H^3_{rad} \times H^2_{rad}, \ \text{supp} (f_j) \subset B(0,R) \}.$$  Define
$$\Sigma_0 = \{ (f_1,f_2) \in X_R : \langle k(\cdot,1) f_1 + f_2, g(\cdot,1) \rangle_{L^2} = 0 \},$$
where $-k^2$ is the unique negative eigenvalue of the linearized
Hamiltonian and $g$ is the associated ground-state eigenfunction.
There exists $\delta > 0$ such that if $(f_1,f_2) \in B_\delta(0)
\subset \Sigma_0$, then there is a real number $h(f_1,f_2)$ such
that
\begin{eqnarray*}
|h(f_1,f_2)| & \lesssim & \|(f_1,f_2)\|^2_{X_R}, \\
| h(f_1,f_2) - h(\tilde{f_1},\tilde{f_2}) | & \lesssim & \delta \| (f_1,f_2) - (\tilde{f_1},\tilde{f_2}) \|_{X_R} \\
 \text{for} \ (f_1,f_2),(\tilde{f_1},\tilde{f_2}) & \in & B_\delta(0)
\end{eqnarray*}
and
\begin{eqnarray*}
\left\{ \begin{array}{l}
\Box \psi - \psi^5 = 0 \\
\left( \begin{array}{l}
\psi(\cdot,0) \\
\partial_t \psi(\cdot,0)
\end{array} \right) = \left( \begin{array}{l}
\phi(\cdot,1) \\
0
\end{array} \right) + \left( \begin{array}{l}
f_1(\cdot) + h(f_1,f_2) g( \cdot,1) \\
f_2(\cdot)
\end{array} \right)
\end{array} \right.
\end{eqnarray*}
has a global solution. Moreover, this global solution admits the
representation $\psi(\cdot,t) = \phi( \cdot, a(\infty) ) +
v(\cdot,t)$ where
\begin{eqnarray*}
\| v( \cdot,t) \|_{\infty} & \lesssim & \delta \langle t \rangle^{-1}, \\
\| (v,\partial_t v) - (\tilde{v},\partial_t \tilde{v}) \|_{\dot{H}^1
\times L^2} & \to & 0
\end{eqnarray*}
as $t \to \infty$  where $\Box\tilde v=0$ and $(\tilde{v},\partial_t
\tilde{v}) \in \dot{H}^1 \times L^2$.
\end{thm}

The radial assumption is more than a mere convenience. It assures
that the bulk term remains at the origin. In other words, in case of
nonradial perturbations we would need to allow $\phi$ to move away
from the origin, i.e., $\phi(X(t,x),a(t))$ where $X(x,t)$ preserves
the Lorentz invariance. In the same spirit, recall that the
eigenspace at zero is generated by the functions $\nabla\phi$ (an
expression of translation invariance), whereas the resonance at zero
results from the dilation invariance. The latter is the only allowed
symmetry in the radial case.

The map
\[ B_\delta(0)\to H^3_{\rm rad}\times H^2_{\rm rad}:\; (f_1,f_2)\mapsto (f_1(\cdot) + h(f_1,f_2) g( \cdot,1),f_2) \]
traces out a (small) manifold $\Sigma$. By our estimates, $\Sigma$
can be written as a Lipschitz graph and  $\Sigma_0$ is the tangent
plane to the desired manifold $\Sigma$ (because of the quadratic
bound on~$h$).

For a general Schr\"odinger operator $H=-\Delta+V$ in $\R^3$ the
spectrum changes under the  scaling $V\to \lambda^2 V(\lambda x)$ by
multiplication by~$\lambda^2$. By inspection,  $\phi(\lambda x,a) =
\lambda^{-\frac{1}{2}} \phi(x,\lambda^2 a)$. Hence, our linearized
operators \[ H(a) = -\Delta - 5\phi^4(\cdot,a) \] have the property
that the spectrum depends  on~$a$ via scaling. In particular, this
is true of  $k$  (and of course $g(\cdot,a)$ depends on~$a$ in a
similar fashion).

To prove the theorem, we make the ansatz  $\psi(t,x) = \phi(x,a(t))
+ u(t,x)$ with some path that satisfies $a(0)=1$ (we will in fact
also need that $\dot{a} (0) = 0$ as well as $| \dot{a} (t) |
\lesssim  \delta \langle t \rangle^{-2}$). Then we plug this ansatz
into our NLW which yields
\begin{eqnarray*}
\partial_{tt} u - \Delta u - 5 \phi^4( \cdot,a(t)) u(t) =
-\partial_{tt} \phi(\cdot,a(t)) + N(\phi(\cdot,a(t)),u)
\end{eqnarray*}
with the nonlinearity \[ N(u,\phi)= 10u^2\phi^3 + 10 u^3\phi^2 + 5
u^4\phi + u^5\] Let us define $V(a) = -5 \phi^4(\cdot,a)$. We want
to eliminate the time dependence in $V$, so we look at $
V(a(\infty))$ and define
\begin{eqnarray*}
\Box_{V(a(\infty))} u = \partial_{tt} u - \Delta u - 5
\phi(\cdot,a(\infty)) u.
\end{eqnarray*}
Thus, the linearized equation becomes
\[ \Box_{V(a(\infty))} u = -\partial_{tt} \phi(\cdot,a(t)) + N(\phi(\cdot,a(t)),u)
+ (V(\cdot,a(\infty))-V(\cdot,a(t)))u =: (*)\]

Before we say more on how to analyze this equation, let us look at
the Schr\"odinger operator $H(a) = -\Delta + V(\cdot,a)$.  The
following relations were already noted above:
\begin{eqnarray*}
H(a) \partial_a \phi & = & 0 \ \text{from dilation invariance} \\
H(a) \nabla \phi & = & 0 \ \text{from translation invariance}.
\end{eqnarray*}
Since $\partial_a \phi \not\in L^2$  we have a resonance at the
point $0$ of the spectrum.

Let us discuss the spectrum of the above operator (on the radial
subspace). First of all, since $V(a)$  decays like $|x|^{-4}$,  a
result of Agmon and Kato implies that there is no singular
continuous spectrum as well as no imbedded eigenvalues in the
continuous spectrum. We also know that the absolutely continuous
spectrum is the same as that of $-\Delta$, namely $[0,\infty)$.
Hence, we can decompose $L^2_{\rm rad} = L^2_{pp} \oplus L^2_{ac}$.
The pure-point part $L^2_{pp}$ is a  one-dimensional space
consisting of $\text{span} \{ g\}$.  Since $g$ is a ground-state
solution, it decays exponentially. This decomposition is blind to
the resonance. In particular, there is no $L^2$ projection
associated with the resonance function. We shall say more about the
characterization of resonances in terms of the perturbed resolvent
in the next lecture. The basic principle is as follows: Look at
$(H(a) - z^2)^{-1}$. As $z \to 0$ in $\Im z>0$, this can either
approach a bounded or an unbounded operator in a suitable weighted
$L^2$ topology. Due to the work of Jensen and Kato from 1978, we
know that this operator can be Laurent expanded in the form
\begin{eqnarray*}
(H(a) - z^2)^{-1}=\frac{c_{-2}}{z^2} + \frac{c_{-1}}{z} + O(1).
\end{eqnarray*}
Clearly, $-c_2$ must equal the orthogonal projection onto the zero
energy eigenspace, whereas $c_{-1}$ has contributions from both the
resonance and the eigenfunctions at zero.  For now, let us not
concern ourselves with respect to which topology this is valid. Let
us rather note that in our application we not only eliminate
$c_{-2}$ but also that $c_{-1}$ is solely due to the zero energy
resonance (recall the restriction to the radial subspace).

We now return to a more detailed description of the main ansatz.
Rewrite the linearized equation as a Hamiltonian system:
\begin{eqnarray*}
\partial_t U = J \mathcal{H}_\infty U + \left( \begin{array}{l}
0 \\
(*)
\end{array} \right),
\end{eqnarray*}
where
\begin{eqnarray*}
U = \left( \begin{array}{l}
u \\
\partial_t u
\end{array} \right),
J = \left( \begin{array}{cc}
0 & 1 \\
-1 & 0
\end{array} \right),
\end{eqnarray*}
\begin{eqnarray*}
\mathcal{H}_{\infty} = \left( \begin{array}{cc}
H(a(\infty)) & 0 \\
0 & 1
\end{array} \right).
\end{eqnarray*}

It is of course necessary to understand the propagator
$e^{tJ\mathcal{H}_{\infty}}$, and therefore also the spectrum of the
operator $J\mathcal{H}_\infty$. To this end, note that if
\begin{eqnarray*}
J \mathcal{H}_{\infty} \left( \begin{array}{l}
f_1 \\
f_2
\end{array} \right) & = & E \left( \begin{array}{l}
f_1 \\
f_2
\end{array} \right), \\
\text{ then\ \ \  }H(a(\infty)) f_1 & = & - E^2 f_1.
\end{eqnarray*}
The conclusion here is that $\text{spec} (J \mathcal{H}_{\infty}) =
\sqrt{- \text{spec}H(a(\infty))} = i \reals \cup \{ \pm k \}$.
Moreover, both $\pm k$ are simple eigenvalues with
eigenfunctions~$G_{\pm}$.

Define $P_{-}$ and $P_{+}$ to be  the Riesz projections associated
with the eigenvalues $-k$ and $k$, respectively.  This can be done
in the standard way by taking contour integrals of the resolvent.
More importantly, they turn out to be very explicit rank one
operators with kernels defined by the ground state function and
eigenvalue.  Decompose the full solution $U(\cdot,t)$ as
\begin{eqnarray*}
U(\cdot,t) = n_{+} (t) G_{+} + n_{-} (t) G_{-} + \tilde{U}(\cdot,t),
\end{eqnarray*}
where $(I - P_{+} - P_{-}) U = \tilde{U}$.  Then, we get the
following hyperbolic system of ODEs in the plane  (with
$k=k(\cdot,a(\infty))$):
\begin{eqnarray*}
\left( \begin{array}{l}
\dot n_{+}(t) \\
\dot n_{-}(t)
\end{array} \right) + \left( \begin{array}{cc}
-k & 0 \\
0 & k
\end{array} \right) \left( \begin{array}{l}
n_{+}(t) \\
n_{-}(t)
\end{array} \right) = \left( \begin{array}{l}
F_{+}(t) \\
F_{-}(t)
\end{array} \right) \in L^\infty_t (0,\infty).
\end{eqnarray*}
Now let us find the unique line in the plane with the property that
initial data chosen on it lead to bounded solutions:
\begin{eqnarray}
n_{+} (t) & = & e^{tk} n_{+} (0) + \int_0^t e^{k(t-s)} F_{+} (s)\, dx \nonumber \\
e^{-tk} n_{+} (t) \to 0 & = & n_{+} (0) + \int_0^\infty e^{-ks} F_{+} (s)\, ds \label{eq:stab}\\
\Rightarrow n_{+} (t) & = & -\int_t^\infty e^{k(t-s)} F_{+} (s)\,
ds. \nonumber
\end{eqnarray}
It follows from the final relation that $n_+(t)$ inherits the decay
of $F_+(t)$ (e.g.~any decay like $ t^{-\beta}$ with $\beta>0$)
provided $n_+(0)$ is chosen as in~\eqref{eq:stab}. Thus,
\eqref{eq:stab} is precisely the stability condition that we use to
define the (unique) correction to the initial data which leads to
bounded $U$ at each step of the Banach iteration. Summing up all
these corrections then yields the function~$h(f_1,f_2)$. Next, by
the explicit form of the projections $P_{\pm}$ one checks that
\begin{eqnarray*}
\tilde{U} = \left( \begin{array}{l}
\tilde{u} \\
\partial_t \tilde{u}
\end{array} \right)
\end{eqnarray*}
where $\tilde{u} = P^{\bot}_{g(\cdot,a(\infty))} u$. Hence, we can
go back from the Hamiltonian system to the wave equation
\begin{eqnarray*}
\Box_{V(\cdot,a(\infty))} \tilde{u} & = & P^{\bot}_{g(\cdot,a(\infty))} (*), \\
\tilde{u} (0) & = & P^{\bot}_{g(\cdot,a(\infty))} (f_1 + h(f_1,f_2) g(\cdot,1)), \\
\partial_t \tilde{u} (0) & = & P^{\bot}_{g(\cdot,a(\infty))} f_2
\end{eqnarray*}
where $(*)$ stands for the right-hand side of the linearized
equation, see above. As a starting point for bounding the solutions
of this wave-equation,  we show that for some constant $c_0\ne0$
\begin{equation}
\label{eq:sin_evol} \frac{\sin (t
\sqrt{H(a(\infty))})}{\sqrt{H(a(\infty))}}
P^{\bot}_{g(\cdot,a(\infty))} f =  c_0 \partial_a
\phi(\cdot,a(\infty)) \otimes \partial_a \phi(\cdot,a(\infty)) f +
\mathcal{S}_\infty (t) f,
\end{equation}
where  the first operator is due to the resonance (and therefore
must be rank one) and $\| \mathcal{S}_{\infty} (t) f\|_\infty
\lesssim \langle t \rangle^{-1} \| f \|_{W^{1,1}}$. If we  choose
$f$ to be a Schwartz function, say, then the first operator here is
well-defined (despite the fact that $\partial_a \phi \not\in L^2$)
and it is given explicitly as  \[c_0\,
\partial_a \phi(x,a) \int
\partial_a \phi(y,a) f(y)\, dy.\]  Let us give a naive reason for the constant term
which the resonance produces in the sine-evolution. First, the
sine-evolution operator acting on $\nabla \phi$ (a zero
eigenfunction) would give a growth rate of $t$ (just Taylor-expand).
The dispersive part decays like $t^{-1}$. Naively, the resonance
should be half-way in between and grow something like $t^0$.  This
is actually the case. In our proof, we show also that the cosine
evolution operator is dispersive.

\bigskip
We continue our discussion of the stable manifold theorem for NLW.
As mentioned previously, we make the ansatz
\begin{equation}
\label{eq:aufteil} \psi(x,t)=\phi(x,a(t)) + u(x,t) \end{equation}
The theorem is proved by solving for $a$ and $u$ globally in time in
such a way that $u$ disperses and $a(t)$ converges to a limit which
is $\delta$-close to $a(0)=1$. Note, however, that we cannot prove
energy scattering for $u$ to a free wave, which explains why we
stated the scattering for a different decomposition of $u$. To
understand this better, one checks from the Duhamel
formula\footnote{This was pointed out by Herbert Koch} for the
nonlinear solution $\psi$ that for any $t_2>t_1>0$
\[ \psi(\cdot,t_2)-\psi(\cdot,t_1) \in L^2 \]
Since \[ \psi(\cdot,t_2)-\psi(\cdot,t_1) =
\phi(\cdot,t_2)-\phi(\cdot,t_1) + u(\cdot,t_2)- u(\cdot,t_1)
\]
and $\phi(\cdot,t_2)-\phi(\cdot,t_1) \not\in L^2$  we see that
\[ u(\cdot,t_2)- u(\cdot,t_1) \not\in L^2 \]
On the other hand, if $u$ were to scatter to a free wave, then
$\partial_t u\in L^2$ and consequently $u(\cdot,t_2)- u(\cdot,t_1)
\in L^2$. The truth is that $u$ scatters to a free energy solution
up to a small, decaying piece which has the same behavior at spatial
infinity as $\partial_a \phi(\cdot,a)$. In fact, our proof shows
that
\[ (u,\partial_t u) = (\tilde u,\partial_t \tilde u) +
(0,-\dot{a} (t)
\partial_a \phi ( \cdot, a(t))) + o_{\dot H^1\times L^2}(1)\]
as $t\to\infty$ where $(\tilde u,\partial_t u)$ is a free wave with
$(\tilde u,\partial_t u)(0)\in \dot H^1\times L^2$. Consequently, if
we set
\[ v(\cdot,t) = u(\cdot,t) + \phi(\cdot,a(t))-\phi(\cdot,a(\infty))
\]
then $v$ scatters and
\[ \psi(\cdot,t) = \phi(\cdot,a(\infty)) + v(\cdot,t) \]
as claimed in our theorem. The ansatz~\eqref{eq:aufteil} modifies
the profile of $\phi$ at spatial infinity instantaneously. While
this is convenient technically, it may be less so physically, since
perturbations should propagate at finite speed. The fact that we
modify $\phi$ outside of the light-cone of course forces us to
absorb the change in that region by means of~$u$. This is exactly
why $u$ does not scatter. It is natural to ask for an ansatz which
does modify the profile of $\phi$ only inside the light-cone. Thus,
the question is as follows: can we write the stable solutions as
\[ \psi (\cdot,t) = \phi(\cdot,a(t)) \chi_1 (\cdot,t) + \phi(\cdot,1) \chi_2 (\cdot,t) + u(\cdot,t)\]
 where $\chi_1$ is a cut-off inside the light cone and $\chi_2$ is a cut-off outside the light-cone, and $u$ disperses and scatters
 in the energy space? It remains to be seen if our proof can
 accommodate this ansatz.

One of our goals for this lecture is to explain how to obtain the
ODE for $a(t)$. This hinges on the linear estimates that we derived
for the wave equation with a potential. During the past twelve years
there has been much activity on  dispersive and Strichartz estimates
for the wave equation with a potential. Specifically, we would like
to mention the work by
 Beals, Strauss, Cuccagna, Georgiev, d'Ancona, Pierfelice,  Visciglia, and Yajima.
However, these authors either assume that the potential is positive,
or small, or that zero energy is neither an eigenvalue nor a
resonance. In addition, varying degrees of smoothness and decay are
assumed. Since we are dealing with a  negative potential that leads
to a zero energy resonance, none of these results apply to our
problem.

However, for the Schr\"odinger equation with a potential there has
been recent work of a similar flavor to what is needed here. Let us
mention the following result from~2004 (there is a matrix version of
this result from 2005 which is relevant to the linearized NLS, see
above).

\begin{thm}
(Erdogan-Schlag) \\
For any Hamiltonian $H = -\Delta + V$ with a real-valued potential
$V(x)$ decaying like $\langle x\rangle^{-\beta}$, $\beta>12$,  we
have
\begin{eqnarray*}
\| (e^{itH} P_c - t^{-\frac{1}{2}} F_t) f \|_{\infty} \lesssim
t^{-\frac{3}{2}} \| f \|_1,
\end{eqnarray*}
where $F_t$ is an  operator for which $\sup_t \| F_t \|_{L^1 \to
L^\infty} < \infty$.
\end{thm}

Furthermore, if zero is not an eigenvalue but only a resonance, then
$F_t$ is rank one, and one can give an explicit expression for
$F_t$. Independently, Yajima proved a similar result which required
less decay of $V$ and he gave an explicit expression for $F_t$ in
all cases.

As mentioned before, a unifying approach to resonances and/or
eigenvalues at zero energy is given by the Laurent expansion of the
perturbed resolvent around zero. To this end,  we again recall the
Jensen, Kato formula
\begin{eqnarray*}
(-\Delta + V - z^2)^{-1} = \frac{c_{-2}}{z^2} + \frac{c_{-1}}{z} +
c_0 + c_1 z + ...,
\end{eqnarray*}
where we can incorporate the non-singular terms into $O(1)$ relative
to a suitable weighted $L^2$ operator norm. To understand the
relevance of this expansion, let us consider the free case $V=0$ in
dimensions one and three.

\begin{example}
Let $d = 1$ and $V=0$.  We then have, for all $\Im z>0$,
\begin{eqnarray*}
(-\frac{d^2}{dx^2} - z^2 )^{-1} (x,y) & = & \frac{e^{iz|x-y|}}{2iz} \\
& = & \frac{1}{2iz} + \frac{e^{iz|x-y|}-1}{2iz} \\
& = & \frac{c_{-1}}{z} + O(1).
\end{eqnarray*}
Hence, $c_{-1} = \frac{1}{2i} 1 \otimes 1$.  This tells us that $1$
is a resonance function. Indeed, it is clear that $y(x) = 1$ is a
solution of $- y'' (x) = 0$. More generally, we say that zero is a
resonance of $-\frac{d^2}{dx^2}+V$ provided there is a bounded,
nonzero, solution of
\[ (-\frac{d^2}{dx^2}+V)f=0\]
The reason for $L^\infty$ lies with the mapping properties of the
free resolvent. Also, if $V$ decays sufficiently fast, then zero
energy is never an eigenvalue (this is specific to $d=1$).
\end{example}

\begin{example}
For $d = 3$ and $V = 0$, we have for $\Im z>0$
\begin{eqnarray*}
(-\Delta - z^2)^{-1} (x,y) = \frac{e^{i z |x-y|}}{4 \pi |x-y|}.
\end{eqnarray*}
Hence, we see that there is no singular term in $z$ and hence no
resonance at $0$. As in the case $d=1$ there is a characterization
of resonances as solutions of $Hf=0$. Indeed, we say that zero is a
resonance iff  there is $f \in L^{2,-\sigma} \setminus L^2$ for
$\sigma > \frac{1}{2}$ $(-\Delta + V) f = 0$.
 In particular, $f(x) = 1$ is not a resonance in $d=3$. To understand this definition better,
 let us  write $(-\Delta + V) f = 0$ equivalently as $f + (-\Delta)^{-1} Vf = 0$.
 Explicitly,
\begin{eqnarray*}
f(x)=-\frac{1}{4\pi} \int_{\R^3} \frac{V(y)}{|x-y|} f(y) dy =
-\frac{1}{4\pi} \int_{\R^3} V(y) f(y)\, dy + O(\langle x
\rangle^{-2}).
\end{eqnarray*}
This shows us two things: \begin{itemize} \item Since
$|x|^{-1}:L^{2,\alpha_1} \to L^{2,-\alpha_2}(\R^3)$ as convolution
operator provided $\alpha_1,\alpha_2>1/2$ and $\alpha_1+\alpha_2>1$,
we see that necessarily
\[ f\in L^{2,-\sigma} (\R^3)\quad \forall \; \sigma>\frac12\]
at least if $|V(x)|\lesssim \langle x\rangle^{-2-\epsilon}$. This
explains the definition of a resonance function.

\item In fact, in that case $f$ decays like $|x|^{-1}$ at infinity
iff $\int Vf(y)\, dy\ne0$ (the resonant case), whereas $f$ is an
eigenfunction decaying like $|x|^{-2}$ at infinity iff $\int Vf(y)\,
dy=0$.

\end{itemize}

\end{example}

It is of course important to understand how to obtain Laurent
expansions of the perturbed resolvent. This was done by Jensen and
Kato around 1978, but a more transparent approach was recently found
by Jensen and Nenciu. Let us state the main lemma from their paper.

\begin{lem}
(Jensen-Nenciu - 2001) \\
Let $H$ be a Hilbert space, and suppose that bounded operators
\[A(z) = A_0 + z A_1 (z):H \to H\] are given for all $z \in F \subset
\cx \setminus \{ 0 \}$. In addition, assume that $A_1(z)$ is
uniformly bounded in $z$, and that $A_0^* = A_0$. Suppose that zero
is an isolated point of the spectrum of $A_0$ and let
\[S = \text{Proj} ( \text{Ker} A_0 ) \] be of finite, positive rank.
Then $(A_0 + S)^{-1}$ exists, and thus $(A(z) + S)^{-1}$ also exists
for small $z$. Define \[B(z) = \frac{1}{z} (S - S(A(z) +S)^{-1} S)\]
Then $A(z)$ is invertible iff $B(z)$  is invertible on $SL^2$, and
in that case we have
\begin{eqnarray*}
A(z)^{-1} = (A(z)+S)^{-1} + \frac{1}{z}(A(z)+S)^{-1} S B(z)^{-1}
S(A(z) + S)^{-1}.
\end{eqnarray*}
\end{lem}

Note that $B(z)$ is uniformly bounded in the operator norm as
$z\to0$. Indeed,
\[ S-S(A_0+S)^{-1}S=0 \]
due to the self-adjointness of $A_0$.

 The connection between this lemma and resolvent expansions
is furnished by the symmetric resolvent identity:  Let $V = v^2 U$
were $U = \text{sgn} V$. Then
\begin{eqnarray*}
R_V (z) := ( -\Delta +V - z^2)^{-1} = R_0(z) - R_0 (z) v (U + v
R_0(z) v)^{-1} v R_0 (z)
\end{eqnarray*}
for all $\Im z>0$. Let us set
\[ A(z) = U + v
R_0(z) v = U +  v R_0(0) v + z \frac{v (R_0(z)-R_0(0)) v}{z} =: A_0
+ zA_1(z)
\]
If $A_0$ is invertible, then the Laurent expansion of $R_V(z)$
around $z=0$ has no singular terms. The converse is also true. Thus
zero energy is neither an eigenvalue nor a resonance iff $A_0$ is
invertible. Using Fredholm theory (and the compactness of $v R_0(0)
v$ on $L^2$), we conclude that $A_0$ is invertible iff there is no
nonzero $L^2$ solution of
\[ (U+v(-\Delta)^{-1}v)f  =  0 \]
This in turn is equivalent to the existence of certain solutions $h$
to
\[ (-\Delta+V)h=0 \]
Indeed, formally set $h=(-\Delta)^{-1}vf$.

In our application to NLW, $A_0$ is therefore not invertible. We
apply the Jensen, Nenciu lemma and check that the $B(z)$ is
invertible. Since $S$ is rank-one in this case, this invertibility
is equivalent with the nonvanishing of a number. Not surprisingly,
this number exactly turns out to be $\int V(y)\partial_a \phi(y,a)\,
dy\ne0$, see above. Hence, the Jensen, Nenciu lemma allows us to
perform the Laurent expansion with relatively explicit coefficients.
This latter feature is most relevant for understanding the evolution
operators of the perturbed wave equation. Indeed, by the spectral
calculus,
\begin{eqnarray}
\frac{\sin (t \sqrt{H})}{\sqrt{H}} P_c = \frac{1}{i \pi}
\int^{\infty}_{-\infty} \frac{\sin (t \lambda)}{\lambda} R_V
(\lambda) \lambda\,d\lambda.\label{eq:spec_calc}
\end{eqnarray}
In order to derive dispersive estimates we follow the usual path of
 splitting the integration here into $|\lambda|
> \epsilon$ and $|\lambda| < \epsilon$ for some $\epsilon$ small enough.
It is common knowledge that dispersive estimates for perturbed
operators belong to the realm of a {\em zero energy theory}. The
point here is that energies different from zero cannot destroy the
free rate of decay. Only zero energy can do that. Nevertheless, the
regime $|\lambda|>\epsilon$ is non-trivial in particular with
respect to the amount of decay of $V$ required to make this work.
Recall that we have exactly $|x|^{-4}$ decay in NLW. To appreciate
the difficulty of lowering the decay requirements on $V$ see
Goldberg's paper~\cite{Gold}.
 At any rate,
in this regime we work with the (finite) Born expansion
\begin{eqnarray*}
R_V = R_0- R_0 V R_0 + R_0 V R_0 V R_0 \mp ... \pm R_0 V  \cdots V
R_V V\cdots V R_0
\end{eqnarray*}
Plugging the terms here not involving $R_V$ into the right-hand side
of~\eqref{eq:spec_calc} leads to completely explicit expressions for
the kernels due to the explicit nature of $R_0$ in $\R^3$. It is
therefore possible, with some work, to obtain the desired dispersive
bounds on those kernels. For the remainder term, which still
involves $R_V$, we rely on the so-called limiting absorption
principle due to Agmon. This refers to the fact that the (perturbed)
resolvent remains bounded on weighted $L^2$ spaces. For further
details on not small energies, we refer the reader to our paper.

For energies $|\lambda|<\epsilon$, we expand the resolvent $R_V$ as
in the Jensen, Nenciu lemma. The singular term ultimately leads to
the non-decaying part
\[ c_0\, \partial_a\phi(\cdot,a)\otimes \partial_a \phi(\cdot,a) \]
whereas we show in our paper that the nonsingular terms yield a
dispersive operator $\mathcal{S}(t)$, see \eqref{eq:sin_evol} above.
Dispersion here refers to bounds in $L^\infty$. It is an open
question whether or not $\mathcal{S}(t)$ satisfies Strichartz
estimates. For this see the author's work on Littlewood-Paley theory
and the distorted Fourier transform.

To conclude these lectures, let us turn to the question of how to
determine the ODE for $a(t)$. To this end let us recall that we
write the solution vector in the form
\begin{eqnarray*}
\left( \begin{array}{l}
u(\cdot,t) \\
\partial_t u(\cdot,t)
\end{array} \right) = n_{+}(t) G_{+} + n_{-}(t) G_{-} + \tilde{u}(\cdot,t)
\end{eqnarray*}
As explained before, generally speaking $n_+$ is an exponentially
unstable mode. It needs to be stabilized by means of the stability
condition~\eqref{eq:stab}. This is the origin of the codimension one
condition. On the other hand, $n_-$ is exponentially decaying.
Finally, the remainder $\tilde u$ satisfies the equation
\begin{eqnarray}
&& \partial_{tt} \tilde{u} + H(a(\infty)) \tilde{u} = \nonumber\\
&&P_{g_\infty}^{\bot} [-\partial_{tt} \phi(\cdot,a(t)) +
(V(\cdot,a(\infty))-V(\cdot,a(t))u(\cdot,t)+
N(u(\cdot,t),\phi(\cdot,a(t)) ] \label{eq:utilrhs}
\end{eqnarray}
It is our goal to show that $\tilde u$ (and therefore also $u$)
disperses, but this is non-obvious because of the non-dispersive
term in~\eqref{eq:sin_evol}. This is precisely where $a(t)$ comes
in, the point being that we will collect all non-decaying terms and
observe that they are all multiples of $\psi:=\partial_a
\phi(\cdot,a(\infty))$. To give a flavor of this, let us consider
the contribution of the first term in~\eqref{eq:utilrhs} to the
Duhamel formula. It is (writing $H$ instead of $H(a(\infty))$ for
simplicity)
\[ -\int_0^t \frac{\sin((t-s)\sqrt{H})}{\sqrt{H}}P_{g_\infty}^\perp
\partial_{s} \dot{a}(s) \partial_a \phi(\cdot,a(s))\, ds \]
We could use the representation \eqref{eq:sin_evol} here. However,
this is very misguided since it is impossible to let
\[ c_0\, \psi\otimes \psi \]
act on $\partial_a\phi(\cdot,a(s))$ (because of $\partial_a
\phi\not\in L^2(\R^3)$). Instead, we proceed as follows:
\begin{align}
-\int_0^t \frac{\sin((t-s)\sqrt{H})}{\sqrt{H}}P_{g_\infty}^\perp &
\partial_{s} \dot{a}(s) \partial_a \phi(\cdot,a(s))\, ds
= \dot{a}(0)\frac{\sin(t\sqrt{H})}{\sqrt{H}}P_{g_\infty}^\perp  \partial_a \phi(\cdot,a(0))\nonumber \\
& + \int_0^t   \dot{a}(s) \,{\cos((t-s)\sqrt{H})}P_{g_\infty}^\perp
\partial_a \phi(\cdot,a(s))\, ds \label{eq:adot_int}
\end{align}
Now observe that
\[ \frac{\sin(t\sqrt{H})}{\sqrt{H}}P_{g_\infty}^\perp  \psi
 = t\psi,  \qquad \cos(t\sqrt{H})\psi=\psi \]
Using these relations and the fact that (for $a_1,a_2 \in (1/2,2)$)
\[ \partial_a \phi(\cdot,a_1) = (a_2/a_1)^{\frac54} \partial_a \phi(\cdot,a_2)
+ O(|a_1-a_2|\langle x\rangle^{-3})
\]
we reduce $\partial_a \phi(\cdot,a(s))$ to $\psi$ by means of the
scalar factor $(a(\infty)/a(s))^{\frac54}$. This shows that we need
$\dot{a}(0)=0$, and also allows us to pull out $\psi$ from the
integral above with a factor depending on $a(t)$. The remainder
which comes from the $O(\langle x\rangle^{-3})$--term turns out to
be dispersive.
 We refer the reader to our paper for further details.


\begin{thebibliography}{99}

\bibitem[Agm]{Ag} Agmon, S. {\em Spectral properties of Schr\"odinger
operators and scattering theory.} Ann.\ Scuola Norm.\ Sup.\ Pisa
Cl.\ Sci. (4) 2 (1975), no.~2, 151--218.

\bibitem[Bea]{Bea} Beals, M. {\em Optimal $L\sp \infty$ decay for solutions to the wave equation with a potential.}
  Comm.\ Partial Differential Equations  19  (1994),  no.~7-8, 1319--1369.

\bibitem[BeaStr]{BS} Beals, M., Strauss, W. {\em $L\sp p$ estimates for the wave
equation with a potential.}  Comm.\ Partial Differential Equations
18 (1993),  no.~7-8, 1365--1397.

\bibitem[BerCaz]{BC} Berestycki, H., Cazenave, T.
{\em Instabilit\'e des \'etats stationnaires dans les \'equations de
Schr\"odinger et de Klein-Gordon non lin\'eaires.} (French. English
summary) [Instability of stationary states in nonlinear
Schr\"odinger and Klein-Gordon equations] C.\ R.\ Acad.\ Sci.\ Paris
S\'er.~I Math.\ 293 (1981), no.~9, 489--492.

\bibitem[BerLio]{BL} Berestycki, H., Lions, P.\ L.  {\em Existence d'oudes
solitaires daus les problemes nonlineares du type Klein-Gordon.}
C.R. Acad. Sci. 288 (1979), no. 7, 395--398.

\bibitem[BizChmTab]{Biz}  Bizo\'n, P., Chmaj, T.,
Tabor, Z. {\em On blowup for semilinear wave equations with a
focusing nonlinearity.}  Nonlinearity  17  (2004),  no.~6,
2187--2201.

\bibitem[BusPer1]{BP1} Buslaev, V.\ S., Perelman, G.\ S. {\em Scattering for the
nonlinear Schr\"odinger equation: states that are close to a
soliton.} (Russian) Algebra i Analiz  4  (1992),  no.~6, 63--102;
translation in  St.\ Petersburg Math.\ J.~4 (1993),  no.~6,
1111--1142.

\bibitem[BusPer2]{BP2} Buslaev, V.\ S., Perelman, G.\ S. {\em On the stability of solitary waves
for nonlinear Schr\"odinger equations.  Nonlinear evolution
equations,}  75--98, Amer.\ Math.\ Soc.\ Transl.\ Ser.~2, 164,
Amer.\ Math.\ Soc., Providence, RI, 1995.

\bibitem[BusPer3]{BP3} Buslaev, V.\ S., Perelman, G.\ S. {\em Nonlinear scattering: states that are close to a soliton.}
 (Russian)  Zap.\ Nauchn.\ Sem.\ S.-Peterburg.\ Otdel.\ Mat.\ Inst.\ Steklov.\ (POMI)  200  (1992),  Kraev.\
  Zadachi Mat. Fiz. Smezh. Voprosy Teor. Funktsii. 24, 38--50, 70, 187;
  translation in  J.\ Math.\ Sci.\  77  (1995),  no.~3, 3161--3169.

\bibitem[CazLio]{CaL} Cazenave, T., Lions, P.-L. {em Orbital stability
of standing waves for some nonlinear Schr\"odinger equations.} Comm.
Math. Phys. 85 (1982), 549--561.

\bibitem[Cof]{Cof} Coffman, C. V. {\em Uniqueness of positive solutions of
$\Delta u - u+u^{3}=0$ and a variational characterization of other
solutions.} Arch.\ Rat.\ Mech.\ Anal.\ 46 (1972), 81--95.

\bibitem[Cuc]{Cuc} Cuccagna, S.
{\em On the wave equation with a potential.}
 Comm.\ Partial Differential Equations  25  (2000),  no.~7-8, 1549--1565.

\bibitem[D'AnPie]{DAP} D'Ancona, P., Pierfelice, V.  {\em Dispersive estimate for the wave equation with a rough
potential,} to appear in JFA

\bibitem[DemSch]{DemSch} Demanet, L., Schlag, W. {\em Numerical
verification of a gap condition for linearized NLS}, preprint 2005.

\bibitem[ErdSch1]{ES1} Erdo\smash{\u{g}}an, M. B., Schlag, W.
{\em Dispersive estimates for Schr\"{o}dinger operators in the
presence of a resonance and/or an eigenvalue at zero energy in
dimension three: I}, Dynamics of PDE, vol.~1, no.~4 (2004),
359--379.

\bibitem[ErdSch2]{ES2} Erdo\smash{\u{g}}an, M. B., Schlag, W.
{\em Dispersive estimates for Schr\"{o}dinger operators in the
presence of a resonance and/or an eigenvalue at zero energy in
dimension three: II}, preprint 2005.

\bibitem[FroGusJonSig]{FGJS}  Fr\"ohlich, J., Gustafson, S., Jonsson, B.\ L.\ G., Sigal, I.\ M.
{\em Solitary wave dynamics in an external potential.}  Comm.\
Math.\ Phys.~250  (2004),  no.~3, 613--642.

\bibitem[FroTsaYau]{FTY} Fr\"ohlich, J., Tsai, T.\ P., Yau, H.\ T.
{\em On the point-particle (Newtonian) limit of the non-linear
Hartree equation.} Comm.\ Math.\ Phys.~225  (2002),  no.~2,
223--274.

\bibitem[GeoVis]{GV} Georgiev, V., Visciglia, N. {\em Decay estimates for the wave equation with potential.}
  Comm.\ Partial Differential Equations  28  (2003),  no.~7-8, 1325--1369.

\bibitem[Gol]{Gold} Goldberg, M. {\em Dispersive bounds for the
three-dimensional Schr\"odinger equation with almost critical
potential}, preprint 2004, to appear in GAFA.

\bibitem[GolSch]{GS} Goldberg, M., Schlag, W. {\em Dispersive estimates
for Schr\"odinger operators in dimensions one and three.} Comm.\
Math.\ Phys.\ 251 (2004), no.~1, 157--178.


\bibitem[Gri]{Grill} Grillakis, M. {\em
Analysis of the linearization around a critical point of an infinite
dimensional Hamiltonian system.} Comm.\ Pure Appl.\ Math.~41 (1988),
no.~6, 747--774.

\bibitem[GriShaStr1]{GSS1} Grillakis, M., Shatah, J., Strauss, W.
{\em Stability theory of solitary waves in the presence of symmetry.
I.} J.\ Funct.\ Anal.~74  (1987),  no.~1, 160--197.

\bibitem[GriShaStr2]{GSS2} Grillakis, M., Shatah, J., Strauss, W.
{\em Stability theory of solitary waves in the presence of symmetry.
II.} J.\ Funct.\ Anal.~94  (1990), 308--348.


\bibitem[JenKat]{JK} Jensen, A., Kato, T. {\em Spectral properties of
Schr\"odinger operators and time-decay of the wave functions.} Duke
Math.\ J.\ 46  (1979), no. 3, 583--611.

\bibitem[JenNen]{JenNen} Jensen, A., Nenciu, G.
{\em A unified approach to resolvent expansions at thresholds.}
Rev.\ Math.\ Phys.\ 13 (2001), no.~6, 717--754.

\bibitem[KriSch1]{KrSch} Krieger, J., Schlag, W. {\em Stable manifolds for all monic supercritical NLS in one dimension},
preprint 2005.

\bibitem[KriSch2]{KrSch2} Krieger, J., Schlag, W. {\em Non-generic blow-up solutions for the critical focusing NLS in
1-d}, preprint 2005.


\bibitem[KriSch3]{KrSch3} Krieger, J., Schlag, W. {\em On the focusing critical semi-linear wave equation}, preprint 2005.


\bibitem[Kwo]{Kw} Kwong, M. K. {\em Uniqueness of positive solutions
of $\Delta u - u +u^{p}=0$ in $\R^{n}$.} Arch.\ Rat.\ Mech.\ Anal.\
65 (1989), 243--266.

\bibitem[Lev]{Lev} Levine, H. {\em Instability and nonexistence of global solutions to
nonlinear wave equations of the form $Pu\sb{tt}=-Au+{\mathcal
F}(u)$.}  Trans.\ Amer.\ Math.\ Soc.~192  (1974), 1--21.


\bibitem[McLSer]{McS} McLeod, K., Serrin, J. {\em Nonlinear Schr\"odinger
equation. Uniqueness of positive solutions of $\Delta u + f(u) =0$
in $\R^{n}$.} Arch.\ Rat.\ Mech.\ Anal.\ 99 (1987), 115--145.


\bibitem[MerZaa1]{MZ1} Merle, F., Zaag, H. {\em On growth rate near the blowup surface for semilinear wave equations.}
  Int.\ Math.\ Res.\ Not.~2005,  no.~19, 1127--1155.

\bibitem[MerZaa2]{MZ2} Merle, F., Zaag, H. {\em  Determination of the blow-up rate for a critical semilinear wave equation.}
 Math.\ Ann.~331 (2005), no.~2, 395--416.

\bibitem[MerZaa3]{MZ3} Merle, F., Zaag, H. {\em  Determination of the blow-up rate for the semilinear wave equation.}
 Amer.\ J.\ Math.~125 (2003), no.~5, 1147--1164.


\bibitem[Per1]{Pe} Perelman, G. {\em Some Results on the Scattering of
Weakly Interacting Solitons for Nonlinear Schr\"odinger Equations}
in "Spectral theory, microlocal analysis, singular manifolds",
Akad.\ Verlag (1997), 78--137.


\bibitem[Per2]{Pe2} Perelman, G. {\em On the formation of singularities in solutions
of the critical nonlinear Schr\"odinger equation.} Ann.\ Henri
Poincar\'e 2 (2001), no.~4, 605--673.

\bibitem[Per3]{Pe3}  Perelman, G. {\em Asymptotic stability of multi-soliton solutions for nonlinear Schr\"odinger equations.}
  Comm.\ Partial Differential Equations  29  (2004),  no.~7-8, 1051--1095.

\bibitem[Pie]{Pie} Pierfelice, V. {\em Decay estimate for the wave equation with a small
potential.} preprint 2003, to appear.

\bibitem[RodSch]{RodSch} Rodnianski, I.,  Schlag, W. {\em Time decay for solutions of
Schr\"odinger equations with rough and time-dependent potentials.}
Invent.\ Math.~155 (2004), 451--513.


\bibitem[RodSchSof1]{RSS1} Rodnianski, I., Schlag, W., Soffer, A.
{\em Dispersive Analysis of Charge Transfer Models}, Comm.\ Pure
Appl.\ Math.~58 (2005), no.~2, 149--216.

\bibitem[RodSchSof2]{RSS2} Rodnianski, I., Schlag, W., Soffer, A.
{\em Asymptotic stability of $N$-soliton states of NLS}, preprint
2003

\bibitem[Sch1]{Sch} Schlag, W. {\em Dispersive estimates for Schr\"odinger operators in
dimension two}, Comm.\ Math.\ Phys.~257 (2005), 87--117.

\bibitem[Sch2]{Sch2} Schlag, W. {\em Stable manifolds for an orbitally unstable NLS}, preprint 2004.

\bibitem[Sch3]{Sch3} Schlag, W. {\em A remark on Littlewood-Paley theory for the distorted Fourier transform}, preprint 2005.


\bibitem[Sha]{Sh} Shatah, J. {\em Stable standing waves of nonlinear
Klein-Gordon equations.} Comm.\ Math.\ Phys.\ 91 (1983), no.~3,
313--327.

\bibitem[ShaStr]{ShSt} Shatah, J., Strauss, W.
{\em Instability of nonlinear bound states.}
  Comm.\ Math.\ Phys.~100  (1985),  no.~2, 173--190.


\bibitem[SofWei1]{SofWei1} Soffer, A., Weinstein, M.
{\em Multichannel nonlinear scattering for nonintegrable equations.}
Comm.\ Math.\ Phys.~133 (1990), 119--146

\bibitem[SofWei2]{SofWei2} Soffer, A., Weinstein, M.
{\em Multichannel nonlinear scattering, II. The case of anisotropic
potentials and data.} J.\ Diff.\ Eq.~98 (1992), 376--390.


\bibitem[Str1]{Str1} Strauss, W.  {\em Nonlinear wave equations.}
 CBMS Regional Conference Series in Mathematics, 73. AMS, Providence, RI, 1989.


\bibitem[Str2]{Str2} Strauss, W. {\em Existence of solitary waves in
higher dimensions.} Comm.\ Math.\ Phys.\ 55 (1977), 149--162.


\bibitem[SulSul]{SulSul} Sulem, C., Sulem, P.-L.
{\em The nonlinear Schr\"odinger equation. Self-focusing and wave
collapse.}
 Applied Mathematical Sciences, 139. Springer-Verlag, New York, 1999.

\bibitem[Szp]{Szp} Szpak, N. {\em Relaxation to intermediate attractors in nonlinear wave equations},
Theoretical and Mathematical Physics, Vol.\ 127 (2001), no.~3, pp.
817–-826.


\bibitem[TsaYau]{TY} Tsai, T.-P., Yau, H.-T. {\em Stable directions for excited
states of nonlinear Schr\"odinger equations.} Comm.\ Partial
Differential Equations 27 (2002), no.~11-12, 2363--2402.

\bibitem[Wei1]{Wei1} Weinstein, M.\ I.
{\em Modulational stability of ground states of nonlinear
Schr\"odinger equations.}
  SIAM J.\ Math.\ Anal.~16  (1985),  no.~3, 472--491.

\bibitem[Wei2]{Wei2} Weinstein, M.\ I.
{\em Lyapunov stability of ground states of nonlinear dispersive
evolution equations.} Comm.\ Pure Appl.\ Math.~39  (1986),  no.~1,
51--67.


\bibitem[Yaj1]{Y1} Yajima, K. {\em The $W\sp {k,p}$-continuity of wave operators
for Schr\"odinger operators.}  J.\ Math.\ Soc.\ Japan  47  (1995),
no.~3, 551--581.


\bibitem[Yaj2]{Y2} Yajima, K. {\em Dispersive estimate for
Schr\"odinger equations with threshold resonance and eigenvalue},
preprint~2004, to appear in Comm.\ Math.\ Phys.


\end{thebibliography}
\end{document}